\documentclass[a4paper,10pt]{article}

 \usepackage[latin1]{inputenc}
\usepackage{amsfonts}
\usepackage{amsmath}
\usepackage[all]{xy}
\usepackage{amssymb}
\usepackage{amsthm}
\usepackage[american]{babel}
\usepackage{caption}
\usepackage[OT1,T1]{fontenc}
\usepackage{makeidx}
\usepackage{mathptmx}
\usepackage{multicol}
\usepackage{rotating}
\usepackage{subfigure}
\usepackage{upgreek}
\usepackage{xcolor}
\usepackage{lmodern}
\usepackage{graphicx}
\usepackage{pstricks-add}
\usepackage{graphicx}
\usepackage{pst-fun}
\usepackage{psfrag}

\usepackage{indentfirst}

    \usepackage{pst-plot}
    \usepackage{pst-text}
\newgray{grayteste}{.9}

\newcommand{\R}{\mathbb{R}}

\newcommand{\M}{\overline{M}}

\def\lra{\longrightarrow}

\def\g{\tilde{g}}
\def\sph{\mathbb{S}}

\def\r{\mathcal{R}}

\def\D{\mathbb{D}}
\def\a{\mathcal{A}}

\def\CC{\mathbb{C}}

\def\hyp{\mathbb{H}}

\def\cpt{\mathbf{K}}
\def\engy{\mathcal{E}}
\def\sob{\mathcal{W}^{1}}

\newtheorem{theorem}{\bf Theorem}[section]
\newtheorem{assertion}{\bf Theorem}[]
\newtheorem{lemma}[theorem]{\bf Lemma}
\newtheorem{corollary}[theorem]{\bf Corollary}

\newtheorem{definition}[theorem]{\bf Definition}

\theoremstyle{definition}
\newtheorem{remark}[theorem]{\bf Remark}

\newtheorem{ass}[theorem]{\bf Assertion}{\em }

\date{2018-01-31}

\begin{document}

\title{The Gauss map of a complete non flat minimal surface in $\R^3$ with
finite total curvature.}
\author{L. P. Jorge and F. Mercuri}
\date{}

\maketitle


\begin{abstract}
In a remarkable paper published in 1964 Robert Osserman proved that
 the image of the Gauss map of a complete non flat minimal surface 
in $\mathbb R^3$ with finite total curvature omits at most $3$ points. In this 
paper we prove that the Gauss map of such a minimal immersions
 omits at most $2$ points. This is a sharp result since the Gauss
 map of the catenoid omits exactly two points. The 
proof follow from classical global results about  minimal surfaces into 
Euclidean three dimensional space.

\end{abstract}


\section{Introduction}  \hspace{.5cm}
Let $M$ be a complete minimal surface in $\mathbb{R}^3$ with Gauss map $G\colon 
M \to \mathbb{S}^{2}.$ It is  known that $G$ is a holomorphic map from $M$,  
viewed as a Riemann surface, to the Riemann sphere $\mathbb{S}^2$. A classical 
problem in minimal surface theory is to know which results from the complex 
function theory remain true for the Gauss map. For instance, is there a Picard  
type theorem for the Gauss map? L. Nirenberg conjectured, around $1950$, that 
the Gauss map  of a complete minimal surface $M$ of $\mathbb{R}^{3}$ could not 
omit a neighbourhood of a point unless $M$ was  a plane. Nirenberg's 
conjecture, 
a generalization of Bernstein's theorem, can also be viewed as a weak form of 
Picard's theorem. L.  Nirenberg's conjectured  was settled by Robert Osserman 
\cite{kn:O1} in the following theorem.
\begin{theorem}[R. Osserman]\label{general} If $\varphi \colon \! M \to 
\mathbb{R}^{3}$ is a complete non-flat minimal surface  with Gauss map 
$G\colon M \to\mathbb{S}^{2}$ then  $Y=\mathbb{S}^{2}\setminus G(M)$ has zero 
capacity.\ In particular $G(M)$ is dense.
 \end{theorem}
 Osserman's result was greatly improved  by F. Xavier  \cite{kn:X} who showed 
that $\sharp Y\leq 6$, followed by H. Fujimoto \cite{kn:Fu} sharp upper bound  
$\sharp Y\leq4$ fully answering Nirenberg's conjecture\footnote{The Gauss map 
of Scherk's surface omits four points.}. Here $\sharp Y$ is the cardinality of  
$Y$.

Osserman also considered the Gauss map missing point problem for the
class of complete minimal surfaces of finite total curvature.\ For $\varphi$\ 
in this class he showed, following ideas of Huber,  that
\begin{enumerate}
\item[a.] $M$ is conformally equivalent to $\M \setminus E$\ where $\M$\ is a 
compact Riemann surface and $E = \{w_1,\ldots, w_m\} \subset \M$\ is a finite 
set of points.
\item[b.] The Gauss map $G\colon M \lra \mathbb{S}^{2}$ extends to a smooth 
map, denoted by the same symbol, $G\colon\M \lra \mathbb{S}^2$.
\item[c.] The map $G\colon \M \lra \mathbb{S}^2$\ is a branched covering map.
\end{enumerate}
He proved the following.
\begin{theorem}[R. Osserman]\label{ftc} If $\varphi \colon M\to \mathbb{R}^{3} 
$ is a non-flat complete minimal surface  with finite total curvature,  then 
$\sharp Y\leq 3$.
\end{theorem}
The proofs of the two Theorems are based on complex analysis, the Weierstrass 
representation formula, but the two results are of different nature. While the 
first one is a problem in values distribution theory for holomorphic functions 
on the disk,  the second one is topological in nature.\ Let us make clear the 
last assertion.\ In \cite{kn:B-F-M}\ the authors introduce the class of 
surfaces of {\em finite geometric type}. These are immersions $\varphi \colon 
M\to \mathbb{R}^{3}$\ of a surface $M$\ such that $M$\ is complete in the 
induced metric and
\begin{enumerate}
\item $M$\ is diffeomorphic to a compact surface $\overline{M}$\ minus a finite 
set of points, $E = \{w_1, \ldots, w_m\}$,
\item the Gaussian curvature vanishes only at a finite number of points,
\item the Gauss map $G$\ extends to a smooth branched covering, denoted by the 
same symbol, $G \colon \M\to \mathbb{S}^2$.
\end{enumerate}
The points $w_i$, or sometimes punctured neighbourhoods of these points are
called the {\em ends} on $\varphi$.\ If $w \in E$,\ we can consider the
``tangent space at $w$'',\ $T_w\M := G(w)^{\perp}$.\ Jorge and Meeks showed in 
\cite{kn:J-M} that the orthogonal projection of a suitable punctured 
neighborhood of $w$\ on $G(w)^{\perp}$\ is a finite covering over the 
complement of a disk and its order, $I(w)$,\ is called the {\em geometric 
index} 
of $w$.\ They also proved the {\em total curvature formula}\ that relates the 
degree of $G$, the geometric indexes of the end and the Euler characteristic of 
$\M$.\footnote{This formula was proved for the case of minimal surface, by 
Osserman, as an inequality, and used to prove Theorem \ref{ftc}}\ Using this 
formula, together with the Riemann-Hurwitz formula for the branched covering 
$G\colon M \to \mathbb{S}^{2}$\ Barbosa, Fukuoka and Mercuri proved, in 
\cite{kn:B-F-M},the following theorem.
\begin{theorem} The Gauss map of a surface of finite geometric type omits at 
most three points.
\end{theorem}
A natural question is if  the upper bound $\sharp Y\leq 3$ is sharp or it can 
be improved to $\sharp Y\leq 2.$ Let us observe that in all known examples of 
minimal surfaces of $\mathbb{R}^{3}$ with finite total curvature  $\sharp Y\leq 
2$, achieving $\sharp Y =2$ in the Catenoid.  Some partial results in this 
direction are the following:
\begin{itemize}
\item if $\M$ has genus 0 then $\sharp Y \leq 2$ ( \cite {kn:O},
\cite{kn:B-F-M}),
\item if the total curvature is bigger or equal to $-20 \pi$, then $\sharp Y 
\leq 2$ (\cite{kn:F} and \cite{kn:WX}),
\item if the points in $E$\ are regular, in a suitable sense,
then $\sharp Y \leq 2$ (\cite{kn:B-F-M}).
\end{itemize}

In this notes we use basic facts on the topology of compact surfaces and
strong result of minimal surface theory to prove the following result
\begin{assertion}\label{th1}
 If the Gauss map of a complete minimal surface with finite total curvature 
misses $3$ or more points of the sphere then the surface is a plane.
\end{assertion}

We outline the proof of the theorem 1. First we consider the immersion
$X\colon M\to\R^3$ with one of the missing points being the north pole. In 
this position the 
Gauss map can be considered as map into the complex plane $\CC$, after 
composing with the stereographic projection.  It is well know that the 
universal 
covering space of $M,\ \  \phi: \tilde{M} \to M$\ is conformally equivalent to 
either the complex plane $\CC$ or the unit disk\ $\D$. In the first case the 
theorem follow applying the little Picard theorem to the composition of $\phi$\ 
with the Gauss map of the immersion. So we can suppose that $\tilde{M}$\ is 
conformally equivalent to the unit disk.

Now consider the $M$ in position such that the point $(0,\ 0,\ 1)$ is not a 
missing point of Gauss map and choose three points $y_1,\ y_2,\ y_3$ in the set 
$Y=\sph^2\setminus G(M).$
Pulling back the metric of $M$ by $\phi$ we get a complete minimal immersion 
$\varphi\colon\D\to\R^3$ given by $\varphi=X\circ\phi$ with Weierstrass 
representation $\{f,\; g\}$ where $g$ is the stereographic projection of  
Gauss map over $\CC.$  The conformal covering of $G(M)$ 
is given by a modular function
$\psi\colon\D\to\CC\setminus\{y_1,\ y_2,\ y_3\}$.  Then we can lift the Gauss 
map $g\colon\D\to\CC\setminus\{y_1,\ y_2,\ y_3\}$ to a 
conformal map $\tilde g\colon\D\to\D.$ In this way we get a new minimal 
immersion $\tilde\varphi\colon\D\to\R^3$ with Weierstrass representation 
$\{f,\; \tilde g\}.$  Observe that $g$ and $\g$ have the same set of polos each 
one with the same order for both and exactly over the zeros of $f$. The 
order of each zero of $f$ is twice the order of the polos of $g$ or $\g.$  This 
grantee that $\tilde\varphi$ is minimal immersion.  In lemma 
(\ref{lemma3}) we prove that $\tilde\varphi$ is 
complete. Since the image of $\tilde g$ is the unit disk the 
Gauss map of $\tilde\varphi$ is inside one hemisphere of 
$\sph^2.$ The classical stability result of Barbosa-do Carmo \cite{kn:B-doC} 
say that $\tilde\varphi$ is globally stable. Then the image of $\tilde\varphi$ 
is a plane, that is, $\tilde g$ is constant, by do Carmo-Peng \cite{kn:doC-P} 
or Schoen \cite{kn:S}. In particular $g$ is constant and this proves the 
theorem.

Another short conclusion is given by Osserman's theorem (\ref{general}) since 
the image of $\g$ is a half sphere and its complement has not zero capacity 
(see also Remark \ref{remark}).

\begin{remark}  In Fang \cite{kn:F} the author presents an example of complete 
minimal surface immersed into $\R^3$ whose Gauss map misses exactly $3$ points 
but it is a tori with periods and, of course, it
does not have finite total curvature and lemma (\ref{lemma3}) is valid only for 
minimal surfaces of finite geometric type.
\end{remark}
\begin{remark}
Let $M$ be a complete minimal surface with compact boundary at infinite 
$E=\partial_\infty M$ and $M=\M\setminus E$ admitting continuous extension of 
the Gauss map $G\colon\M\to\sph^2.$   
The Fujimoto theorem \cite{kn:F} says that $\sharp Y \leq 4$ where 
$Y=\sph^2\setminus G(M).$ Then $G$ is constant in each connected component of 
$E.$ 
 The proof of completeness of $\tilde\varphi$ uses strongly the fact that 
$G^{-1}(Y)$ has finite number of connected components. It does not work for 
infinite connected components. The same 
construction could be done in any complete minimal surface with Gauss map 
missing $3$ or $4$ points but the completeness of $\tilde\varphi$ depends 
strongly of this fact. The proof is based in the validity of lemma 
(\ref{lemma3}) for this conditions and implies the next results.
There is no complete non flat minimal immersion $X\colon M\to\R^3$ with 
Gauss map $G,$ set of ends $E$ and set of missing points $Y=\sph^2\setminus 
G(M)$ satisfying the following conditions.
 \begin{enumerate}
 \item $\sharp Y=3,$ or $4$,
 \item The Gauss map extends continuously to $\M$,
  \item The set $G^{-1}(Y)$ is compact and has only finite connected 
components.
 \end{enumerate}
 Item 3 implies the degree of $G$ is finite or it has finite total curvature. 
In this particular case the connected components of $E$ are points since $X$ is 
minimal. But this is not the case for general immersions with finite total 
curvature into $\R^3$, for example, harmonic immersions. 
\end{remark}

We would like to thank all people that critically read the manuscript and gave 
helpful suggestions, specially Luciano Mari.

\section{Preliminary results}
Let $M=\M_\mu\setminus E$ be a finite topological type surface with finite
ends $E=\{w_1,\cdots,w_m\}$ and genus $\mu.$
\begin{definition}
An end $N_w,\, w\in E,$ is normal if it is a end and the
boundary curve $\gamma=\partial E_w$ is an embedded closed geodesic of $M.$
\end{definition}
Consider an open set $U\supset E$ given by
\[
U=\bigcup_{j=1}^m U_j,
\]
where $\{U_j\},\, j=1,\cdots,m,$ is a collection of disjoint
ends. Take $U_0=M\setminus U.$ We always consider the boundary curves
$\gamma_j=\partial U_j$ with the orientation induced from $U_0.$
Denote by $\eta_j$ the unit vector field orthogonal to $\gamma_j$ and appoint 
outside $U_0.$
\begin{lemma}\label{lemma1}
Let $M=\M\setminus E$ a finite topological type surface with set of ends 
$E=\{w_1,\cdots,w_m\}.$
There is a collection of normal ends $N_j\ni w_j,\, 1\leq j\leq
m,$ such that
 \begin{enumerate}
\item If $M$ is not an annulus then all the ends $N_j$ are disjoint
\item Two normal ends $N_i$ and $N_j$ have non empty intersection  if
and only if $M$ is an annulus. In that case $\partial N_i=\partial N_j.$
\item The  map $\varphi_j\colon \gamma_j\times [0,\infty)\to
N_j,\; \gamma_j=\partial N_j,$ given by
\[\varphi_j(\theta,t)=\exp_{\gamma_j(\theta)}(t\eta_j(\theta))\]
is a diffeomorphism.
\item If $M$ is not annulus then exist some constant $c>0$ such that if 
$\gamma(t),\;t_1\leq t\leq t_2,$ has $\gamma(t_1)\in N_i$ and $\gamma(t_2)\in 
N_j,\; i\neq j,$ then the length $\ell(\gamma)\geq c.$
 \end{enumerate}

\end{lemma}
\begin{proof} Let $\varphi\colon M\to\R^3$ be a finite geometric type
immersion.Let $A$ be the annulus $1/2\leq|z|\leq2,\, z\in\CC,$ and
$\CC_0=\CC\setminus\{0\}.$ Let $\beta(t),\,0\leq t\leq 1,$ be a
closed curve in $\CC_0$ with winding number $1.$ If $|\beta(t)|\geq 1,\,
\beta(t)=r(t)\exp(\theta(t)),$ then the length $\ell(\beta)$
satisfies
\begin{equation}\label{eq2.1}
\ell(\beta)=\int_0^1\sqrt{|r'|^2+r^2|\theta'|^2}
dt\geq\Big|\int_0^1\theta'dt\Big|\geq2\pi
\end{equation}
If $\beta(t)$ cross the circles $|z|=1$ and $|z|=2$ at the
points $t_0$ and $t_1$ then
\begin{equation}\label{eq2.2}
 \ell(\beta)\geq\int_{t_0}^{t_1}|r'(t)|dt\geq1.
\end{equation}
Hence, any curve $\beta(t)$ with winding number $1$ in $\CC_0$
crossing the circles $|z|=1$ and $|z|=2$ or with $|\beta(t)|\geq1$ has
length $\ell(\beta)\geq 1.$ We know from \cite{kn:J-M} that there is
$R_0>0$ such that
\begin{enumerate}
\item $\varphi^{-1}\big(M\setminus B_{R_0}^{\R^2}(0)\big)$ is a collection of 
ends $U_j\ni w_j,\; 1\leq j\leq m,$ such that $\frac{1}{R_0}\varphi(U_j)$ is a 
multi-graphic over the region $|z|\geq 1/2$ in the plane at infinite,
\item That multigraphic converges with multiplicity the geometric index
$I(w_j)$ in the $C^3$ topology to the annulus $A$ when $R_0\mapsto\infty.$
\end{enumerate}
Since
\begin{equation}\label{eq2.3}
 \frac{1}{R_0}\varphi\big(U_j\big)\mapsto A
\end{equation}
in the $C^3$ topology with multiplicity $I(w_j)$ we can assume that
$R_0$ is bigger enough to have any curve $\alpha(t)$ inside $U_j$ homotopic to 
$\partial U_j$ satisfying
\begin{equation}\label{eq2.4}
\ell(\frac{1}{R_0}\alpha)\geq\ell(\frac{1}{R_0}P(\alpha))\geq
2\pi I(w_j)\geq \pi,
\end{equation}
where $P$ is the orthogonal projection of $\R^3$ over the plane
at infinite $G(w_j)^\perp.$ Take $\alpha\colon I\to M$ a $C^1$ curve with 
extremes at $U_0$ where $I$ is a compact interval. Also we can assume in the 
case of $\alpha(t_1)\in \partial U_j$ and $\alpha(t_0)\in U_0$ that
\begin{equation}\label{eq2.5}
 \ell(\alpha)_{t_0}^{t_1}\geq I(w_j)\geq 1/2.
\end{equation}
The energy of a closed curve $\alpha\colon\sph^1\to M$ is given by
\begin{equation}\label{eq2.6}
 \engy(\alpha)=\int_{\sph^1}|\alpha'|_M^2,
\end{equation}
and Schwartz inequality gives
\begin{equation}\label{eq2.7}
 \ell(\alpha)\leq\sqrt{2\pi}\sqrt{\engy(\alpha)}.
\end{equation}
For a compact set $\cpt\subset M$ can be covered by a finite number of
coordinates and define the norm $\|\alpha\|_{1,\cpt}$ of the Sobolev space
$\sob(\cpt)$. Then we got the following result.
\begin{ass}\label{ass1}{\em
 Given $C>0$ and $\alpha\colon\sph^1\to M$ with $\engy(\alpha)\leq C$ we have 
for all $R>\max\{R_0,\sqrt{2\pi C}\}$ that $\alpha(\sph^1)\subset \cpt_R$ where
\[
 \cpt_R=M\setminus \varphi^{-1}\big(B_R^{\R^3}(0)\big).
\]
Further, there exist some constant $c=c(\cpt_R,C)$ such that
\[
 \|\alpha\|_{1,\cpt_R}\leq c.
\]}
\end{ass}
Let $X_j$ be the class of curves $\alpha\colon\sph^1\to M$ free homotopic to 
$\partial U_j$ and belonging the Sobolev space $\sob(\cpt_M).$ We choose the 
constant $C$ such that
\[
 C>\engy_j=\inf_{\alpha\in X_j}\engy(\alpha).
\]
Fix some $R>0$ given by assertion (\ref{ass1}). Let $X_j(C)$ be the subset of 
$\alpha\in X_j$ with $\engy(\alpha)\leq C.$ Let $\alpha_k\in X_j(C)$ be a 
minimizing sequence with $\lim_k\engy(\alpha_k)=\engy_j.$ By Sobolev immersion 
theorem we can assume that there is $\alpha\in X_j(C)$ such that
$\alpha_k\mapsto\alpha$ in the $C^{0,\nu}$ topology for $0<\nu<1/2.$ We can not 
have $\alpha$ constant since it is not homotopic to a constant. We also have 
that $\alpha_k$ converges weakly in $\sob(M)$ to $\alpha$ and
\[
 \engy(\alpha)\leq \liminf \engy(\alpha_k).
\]
Hence $\engy(\alpha)=\engy_j>0$ and the Euler equation for the energy gives 
that $\alpha$ has the same class of differentiability of $M.$ The standard 
procedure shows that $\alpha$ is one embedded geodesic. Let call this closed 
embedded geodesic $\gamma_j$ and it defines an end $E_j\ni w_j.$

 By the Gauss-Bonnet theorem we can not have that $N_j\cap N_k$ has an 
annulus or a disc as connected component. Then $N_j\cap N_k,\, j\neq k,$ is 
empty or has touching points at the boundary. Since the boundary are geodesics 
then $N_j$ and $N_k,\,j\neq k,$ are not disjoint only if $M$ is an annulus 
proving items $1$ and $2$.

The map defined into item $3$ is over by completeness of $M$ and is one-to-one 
by Gauss-Bonnet theorem. The constant $c$ in item $4$ is the minimal distance 
between the par of closed geodesics $\gamma_j'$s.  

\end{proof}

Let $M=\M\setminus K$ be a connected complete minimal surface into $\R^3$ where 
$K$ is a compact set with finite connected components $K_1,\ \cdots,\ K_m$. 
In similar way we get the existence of normal neighbourhoods $N_j\supset K_j$ 
satisfying the same conclusion of lemma (\ref{lemma1}). 

\begin{corollary}
 If $M=\M\setminus E$ is a connected complete minimal surface with $E=K$ a 
compact set with finite connected components then each connected component 
$K_j$ of $K$ is inside a normal end $N_j$ satisfying 
the items of lemma(\ref{lemma1}). 
\end{corollary}

\section{The basic lemma}\label{secao1}
Let $\hyp=\{z=u+iv\in\CC\,|\,v>0\}$ be the hyperbolic space and
$\hyp^\infty$ be the hyperbolic space with the boundary at infinity.
Let $M=\M_\mu\setminus E$ be a complete minimal surface with compact set of 
ends $E,\ \M$ of genus $\mu\geq0,$ Gauss map $G\colon
M\to\sph^2$ and set of missing points $Y=\sph^2\setminus G(M).$ Since there is 
continuous extension of the Gauss map of $M$ to $\M$ 
we denote both by $G.$  We can consider, after rotation of the immersion if 
necessary, that $y_0=(0,\ 0,\ 1)\in Y.$ Projecting the Gauss
map $G$ into $\CC$ by the stereographic map we get a map $g\colon M\to\CC$
missing at least two points $a$ and $b,\, a\neq b.$ There are two connected 
components of $E$ going into the points  $a$ and $b.$ If we compose $g$ with 
the 
map
\[
 z\mapsto \frac{z-a}{b-a},\qquad z\in\CC,
\]
we get a map from $M$ to $\CC$ denoted yet by $g$ that omit the points $0$ and 
$1$ of the complex plane. By simplicity we may assume that $g$ omits exactly 
$0$ and $1.$ Consider the region $V_0$ of $\hyp^\infty$ given
by
\[
 V_0=\Big\{z\in\hyp^\infty\,\Big|\, 0\leq
u\leq1,\,v\geq\sqrt{1/4-(u-1/2)^2}\;\Big\}.
\]
The boundary of $V_0$ is the union of lines $L_0:u=0,\, v\geq0$
and $L_1:u=1,\, v\geq0$ and the half circle of $C_{1/2}:0\leq u\leq 1,\,
v=\sqrt{1/4-(u-1/2)^2}.$ The classical modular function
$\psi\colon\hyp\to\CC\setminus\{0,\, 1\}$ with basic domain $V_0$
 is done applying $V_0$ into $\hyp^\infty$ by some conformal map $\psi_0\colon 
V_0\to\hyp^\infty$ with $0,\,1$ fixed and after doing all reflections possible 
by the lines and circles at the boundary ever and ever. The first reflection 
over $L_0$ extends
$\psi_0$ to a map $\psi_1\colon V_1\to\CC\setminus\{0,\, 1\}$ where
\[
 V_1=V_0\cup V_0^-,\quad V_0^-=\Big\{ x+iy\in\CC\,\Big|\, -x+iy\in V_0\Big\}
\]
One fundamental region is given by $\Omega_0=V_1\setminus L_{-1}\cup C_{1/2}.$ 
We have that 
$\psi\colon\hyp\to\CC\setminus\{0,\,1\}$ is a covering map. The reflection of 
$V_0$ over the half circle $C_{1/2}$ will give the set of corners 
$\{0,\,1/2,\,1\}. $ This set in the construction of $\psi$ will generate 
sequences dense in all real axis. The boundary $\partial\hyp^\infty$ is just 
$\R\cup\{\infty\}.$ Let $U_0\supset E^\infty= G^{-1}(y_0)$ and $U_1\supset 
E\setminus E^\infty$ disjoints open sets. Let $y_1\in Y$ the possible fourth 
point and choose some $\epsilon_0>0$ such that 
\begin{equation}\label{eq1.1}
 \epsilon_0< \inf\{1/4,\ |y_1|/4,\ |y_1-1|/4,\ (4|y_1|+4)^{-1}\}.
\end{equation}
If 
\begin{equation}\label{eq1.2}
 G(U_1)\subset \bigcup_{w\in Y_1}B_{\epsilon_0}^\CC (w),\qquad Y_1=Y\setminus 
\{y_0\},
\end{equation}
then $g|U_1$ is bounded.  
The horoball
\[
 V(\infty,\epsilon_0)=\big\{ z\in\hyp\;|\;\Im(z)\geq1/\epsilon_0\big\}
\]
with boundary the horosphere $\sigma_{\epsilon_0}=\{(x,\;
1/\epsilon_0)\;|\;x\in\R\}$ is a neighbourhood at $\infty\in\hyp^\infty.$
 The others neighbourhoods at $\infty$ are horoball with vertex over the 
sequence of points reflecting $1/2$ on the construction of $\psi.$

\begin{figure}[!htb]
\centering
\psset{unit=0.3bp}
\begin{pspicture}(-270,-72)(270,256)
\psline(-270,0)(270,0)
\psline(-256, 0)(-256,256)
\psline(256,0)(256,256)
\psline(0,-72)(0,256)

\psline[linecolor=black,linestyle=dashed](-270,170)(270,170)

\psarc(-128,0){128}{0}{180}
\psarc(128,0){128}{0}{180}

\psarc(192,0){64}{0}{180}
\psarc(64,0){64}{0}{180}
\psarc(-64,0){64}{0}{180}
\psarc(-192,0){64}{0}{180}

\psarc(-224,0){32}{0}{180}
\psarc(-160,0){32}{0}{180}
\psarc(-96,0){32}{0}{180}
\psarc(-32,0){32}{0}{180}
\psarc(224,0){32}{0}{180}
\psarc(160,0){32}{0}{180}
\psarc(96,0){32}{0}{180}
\psarc(32,0){32}{0}{180}

\psarc(-240,0){16}{0}{180}
\psarc(-208,0){16}{0}{180}
\psarc(-176,0){16}{0}{180}
\psarc(-144,0){16}{0}{180}
\psarc(-112,0){16}{0}{180}
\psarc(-80,0){16}{0}{180}
\psarc(-48,0){16}{0}{180}
\psarc(-16,0){16}{0}{180}
\psarc(240,0){16}{0}{180}
\psarc(208,0){16}{0}{180}
\psarc(176,0){16}{0}{180}
\psarc(144,0){16}{0}{180}
\psarc(112,0){16}{0}{180}
\psarc(80,0){16}{0}{180}
\psarc(48,0){16}{0}{180}
\psarc(16,0){16}{0}{180}

\pscustom[fillstyle=solid,fillcolor=grayteste]
{
\psline(0,256)(0,0)
\psarcn(128,0){128}{180}{0}
\psline(256,256) 
}
\psline[linecolor=black,linestyle=dashed](0,170)(256,170)

\pscustom[fillstyle=solid,fillcolor=grayteste]
{
\psarcn(64,0){64}{180}{0}
\psarc(96,0){32}{0}{180}
\psarc(32,0){32}{0}{180}
\closepath }

\pscustom[fillstyle=solid,fillcolor=grayteste]
{
\psarcn(192,0){64}{180}{0}
\psarc(224,0){32}{0}{180}
\psarc(160,0){32}{0}{180}
\closepath }

\pscustom[fillstyle=solid,fillcolor=grayteste]
{
\psarcn(-128,0){128}{180}{0}
\psarc(-64,0){64}{0}{180}
\psarc(-192,0){64}{0}{180}
\closepath }

\pscustom[fillstyle=solid,fillcolor=grayteste]
{
\psarcn(-224,0){32}{180}{0}
\psarc(-208,0){16}{0}{180}
\psarc(-240,0){16}{0}{180}
\closepath }

\pscustom[fillstyle=solid,fillcolor=grayteste]
{
\psarcn(-160,0){32}{180}{0}
\psarc(-144,0){16}{0}{180}
\psarc(-176,0){16}{0}{180}
\closepath }

\pscustom[fillstyle=solid,fillcolor=grayteste]
{
\psarcn(-96,0){32}{180}{0}
\psarc(-80,0){16}{0}{180}
\psarc(-112,0){16}{0}{180}
\closepath }

\pscustom[fillstyle=solid,fillcolor=grayteste]
{
\psarcn(-32,0){32}{180}{0}
\psarc(-16,0){16}{0}{180}
\psarc(-48,0){16}{0}{180}
\closepath }

\rput(64,-16){$1$}
\rput(32,-16){$\infty$}
\rput(224,-16){$\infty$}
\rput(192,-16){$0$}
\rput(160,-16){$1$}
\rput(128,264){$\infty$}
\rput(-128,264){$\infty$}
\rput(96,-16){$0$}
\rput(-256,-16){$1$}
\rput(256,-16){$1$}
\rput(128,-16){$\infty$}

\rput(-64,-16){$1$}
\rput(-32,-16){$\infty$}
\rput(-224,-16){$\infty$}
\rput(-192,-16){$0$}
\rput(-160,-16){$1$}

\rput(-96,-16){$0$}

\rput(-8,-16){$0$}
\rput(-128,-16){$\infty$}

\rput(-128,190){$V(\infty,\epsilon_0)$}
\rput(-30,154){$1/\epsilon_0$}

\end{pspicture}
\caption{Fundamental domains with points going to $\{0,\ 1,\ 
\infty\}$}\label{ModDomain01}
\end{figure}

 The arc $\sigma_+=\psi(\sigma_{\epsilon_0}([0,\;1]))$ is embedded inside the 
half space $\Im(z)>0$ and has extremes at real axis $\R_-\colon x<0$ and 
$\R_+\colon x>0.$ We have that $\sigma_-=\psi(\sigma_{\epsilon_0}([-1,0]))$ is 
the reflexion of $\sigma_+$ by the real axis.
  Hence $\sigma=\sigma_-\cup\sigma_+$ is a Jordan curve and define a bounded 
topological disk $\cpt_0$ as the closure of its interior. Consider 
$\r_0=\{(x,y)\;|\;-1\leq x\leq 1,\; y\geq1/\epsilon_0\}.$ It is well know that  
$\psi(z)$ is asymptotic to $\cosh(z)$ for $z\in \r_0$ and $|z|\mapsto\infty.$ 
  But we want to restrict the map $\psi$ in the horizontal strip 
$$\a_0=V(\infty,\epsilon_0)\setminus V(\infty,\epsilon_1),\; 
0<\epsilon_1<\epsilon_0<1/4.$$ We have $\psi(\tau(z))=\psi(z),\ z\in 
\a_0\cap\r_0$ for any deck transform $\tau\in \text{Deck}(\psi).$  
Hence 
\begin{equation}\label{eq1} 
\sup_{z\in\a_0}\frac{|\psi(z)|}{|z|}=\sup_{z\in\a_0\cap\r_0}\frac{|\psi(z)|}{|z|
} <1/c_0,\quad c_0>0.
\end{equation}

\begin{figure}[!htb]
\centering
\psset{unit=0.5bp}
\begin{pspicture}(-190,-72)(210,200)

\psline(-190,-56)(-30,-56)
\psline(-174,-56)(-174,80)
\psline(-46,-56)(-46,80)
\psline(-110,-66)(-110,80)
\psarc(-142,-56){32}{0}{180}
\psarc(-78,-56){32}{0}{180}
\psline(-174,-14)(-46,-14)
\psline[linestyle=dashed](-190,-14)(-30,-14)
\pscustom[fillstyle=solid,fillcolor=grayteste]
{
\psline(-174,80)(-174,-14)(-110,-14)(-110,80)
}
\pscustom[fillstyle=solid,fillcolor=grayteste]
{
\psline(-110,80)(-110,-14)(-46,-14)(-46,80)
}

\psline{->}(-10, 32)(42,32)

\psellipse(140,0)(45,30)
\pscustom[fillstyle=solid,fillcolor=grayteste]
{
\psframe(70,-72)(210,80)
}
\pscustom[fillstyle=solid,fillcolor=white]
{
\psellipse(140,0)(45,30)
}
\psline[linecolor=black,linestyle=dashed](-190,60)(-30,60)

\psframe[linecolor=white](70,-72)(210,80)

\rput(-174,-66){$-\oldstylenums{1}$}
\rput(-100,-66){$\oldstylenums{0}$}
\rput(-46,-66){$\oldstylenums{1}$}
\rput(-25,0){$1/\epsilon_0$}

\psdots*[](140,0)(160,0)
\rput(140,-12){$\oldstylenums{0}$}
\rput(160,-12){$\oldstylenums{1}$}
\rput(120,10){$\mathbf{K}_{\oldstylenums{0}}$}
\rput(160,35){$\sigma$}
\rput(15,22){$\psi$}
\rput(-100,90){$\mathcal{R}_{\oldstylenums{0}}$}
\rput(140,90){$\psi(\mathcal{R}_{\oldstylenums{0}})$}
\psline(185,0)(210,0)
\rput(-200,30){$\a_0$}
\rput(-25,67){$1/\epsilon_1$}

\end{pspicture}
\caption{The region $\a_0.$ where 
$|z|\geq c_0|\psi(z)|$. }\label{ModDomain02}
\end{figure}

We can realize $\M_\mu$ as a regular polygon with
$4\mu$ sides, $\mu\geq1.$ The sides $a_1,\cdots,
a_\mu,$ and $b_1,\cdots,b_\mu,$ are ordered by
$a_jb_ja_j^{-1}b_j^{-1}$ and identified according the
orientation. All the vertex $v_1,\cdots, v_{4\mu}$ represents one point
$x_0\in\M_\mu$ and $\M_\mu$ is homeomorphic to $P_\mu/\sim$ with the vertex of 
$P_\mu$ associated to $x_0.$ Every par $a_j,\,a_j^{-1}$ and $b_j,\,b_j^{-1}$ 
define embedded curves $\alpha_j,\,\beta_j$ crossing only at $x_0.$ The curves 
$\{\alpha_j,\beta_j\, |\, 1\leq j\leq\mu\}$ are the generator of the homotopy 
group $\pi_1(\M_\mu,x_0).$ We identify the set of ends $E$ with a set in the 
interior of $P_\mu.$

From now on we will consider the existence of a complete non flat minimal 
surface $M$ with cardinality of the missing set of points of the 
Gauss map $\sharp Y=3$ or $4.$ Set $Y=\{y_1,\ y_2,\ y_3,\ y_4\}$ for $\sharp 
Y=4$ and for $\sharp Y=3$ the point $y_4$ do not exist. Consider the minimal 
immersion $X\colon M\to\R^3$ in position with $y_0=(0,0,1)\notin Y.$  We set 
\begin{enumerate}
 \item $B_{j\delta}=\text{dist}_{\delta}^{\sph^2}(y_j),\ 0\leq j\leq4$,
 \item $U_{j\delta}=G^{-1}\big(B_{j\delta}\big)$,
 \item $U_\delta=\cup_{j=1}^{4} U_{j\delta}$,  
\item $N_j$ is the finite collection of all normal neighbourhood isolating the 
connected components of $E^j= G^{-1}(y_j)$.
 \item $N=\cup N_j$,
 \item $V_\delta=G^{-1}(B_{0\delta}).$
 \end{enumerate}
 Choose some $\delta_0>0$ such that:
\begin{enumerate}
 \item[(a)] The balls $B_{j\delta_0},\ 0\leq j\leq4$ are disjoints,
\item[(b)] For each $j,\ 1\leq j\leq4,$ the set 
$U_{j\delta_0}\Subset N_j$,
\item[(c)] There is a constant $c_1>0$ such that $\text{dist}^M(\partial 
U_{j\delta_0},\ \partial N_j\big)\geq c_1>0$ for all $j\geq1.$
\end{enumerate}
\begin{remark}
 The item (c) is true since $N$ has finite connected components.
\end{remark}

It is well know that the universal covering of $M$ is conformally equivalent to 
the plane $\CC$ or the disk $\D.$ If there exist a conformal map 
$\phi\colon\CC\to M$ then the map $G\circ\phi\colon\CC\to\CC$ is an entire 
function and by the classical little Picard's theorem it is constant if misses 
$2$ or more points. At this case we need $y_0\in Y$ to get $G\circ\phi$ an 
entire function.

Then the universal covering is conformal to the disk $\D,$ that is,
the conformal covering is a conformal map $\phi\colon \D\to M$ and inducing the 
metric it is a local isometry. Let
$\varphi\colon\D\to\R^3$ be the composition of the minimal immersion $X\colon 
M\to\R^3$ with 
$\phi$ and denote by $f,\;g\colon\D\to\CC$ the Weierstrass representation. The 
map $f\colon\D\to\CC$ is holomorphic and $g\colon\D\to\CC$ is meromorphic.  The 
set of points $W_0=(G\circ\phi)^{-1}(y_0)$ are all polos of $g$ and zeros of 
$f$ and the order of zero of $f$ is twice the order of polos of $g.$ 

Let $h$ be one stereographic map such that $h(W)=\{0,\ 1,\ 
\infty\}$ for $W=\{y_1,\ y_2,\ y_3\}.$ Then $\tilde 
\psi\colon\hyp\to\sph^2\setminus W$ given by $\tilde\psi  =h^{-1}\circ\psi\circ 
h$ is the universal covering. Let $\tilde g$ be the lifting of $g$ by $\tilde 
\psi.$ Since $\tilde\psi$ is local diffeomorphism the map $\tilde g$ has a 
polos at the points $W_0$ with the same order of the polos of $g$. Then the par 
$(f,\ \tilde g)$ is a Weierstrass par of some minimal immersion 

\begin{equation}\label{varphitil}
\tilde\varphi\colon\D\to\R^3.  
\end{equation}

The assumption that $E$ is some compact 
subset of $\M_\mu,$ with finite number of connected components is 
equivalently to say that the total curvature is finite or the degree of $G$ is 
finite.

\begin{lemma}\label{lemma3}
 Let $M$ be a complete non flat minimal surface whose Gauss map extends 
continuously to 
a map $G\colon\M_\mu\to\sph^2,\ M=\M_\mu\setminus E,$ and finite degree. If the 
set 
 of missing points $Y$ has $\sharp Y=3,$ or $4$ then 
$\tilde\varphi\colon\D\to\R^3$ is a complete minimal immersion with image of the 
Gauss map inside a half sphere. In particular such immersion does not exist. 
\end{lemma}

\begin{proof} In the definition of $\tilde{\varphi}$ we had proved not only the 
existence but the absences of branch points in way that it is an immersion. 
We need to prove it is complete. Since the image of the Gauss map of 
$\tilde\varphi$ is inside a half-sphere the facts pointed out in the 
Introduction imply the non existence of such minimal immersion $\varphi.$

Let $\gamma(t),\;0\leq t<\omega_0\leq\infty,$ be a divergent curve in $\D$
parametrized by the arc length of the plane. The length 
$\tilde\ell(\gamma)$ of the immersion $\tilde\varphi$ satisfies
\[
 \tilde\ell(\gamma)=\int_0^{\omega_0}|f|(1+|\tilde
g|^2)|d\gamma|\geq\int_0^{\omega_0}|f||d\gamma|.
\]
Then if $\int_0^{\omega_0}|f||d\gamma|=\infty$ we get 
$\tilde\ell(\gamma)=\infty.$ Here $\omega_0$ is the length of 
$\gamma$ with the Euclidean metric. Then we have 
$\tilde\ell(\gamma)=\infty$ in the following situations.
\begin{enumerate}
 \item[(i)]{\em There is $\delta,\ 0<\delta\leq\delta_0,$ such that 
$\phi(\gamma(t))\in M\setminus V_\delta$ for all $t\geq t_0\geq0.$}.
\item[(ii)]{\em There are increasing sequences $t_j<s_j<t_{j+1}\cdots$ such 
that $\phi(\gamma(t_j))\in U_{\delta_0}$ and $\phi(\gamma(s_j))\notin N$.  }
\item[(iii)] {\em There are increasing sequences $t_j<s_j<t_{j+1}\cdots$ such 
that $\phi(\gamma(t_j))\in V_{\delta'}$ and $\phi(\gamma(s_j))\notin V_\delta$ 
for $0<\delta'<\delta\leq\delta_0.$  }
\end{enumerate}
In item (i) we have $g$ bounded over the image of $\alpha=\phi(\gamma)$ implying
\[
 c_2=\sup_{z\in M\setminus V_\delta}|g(z)|<\infty.
\]
Then
\begin{equation}\label{eq1.0}
 \int_{\gamma}|f||dz|\geq 
\frac{1}{1+c_2^2}\int_{t_0}^{\omega_0}|f|(1+|g|^2)|d\gamma|=\infty,
\end{equation}
and we get $\tilde\ell(\gamma)=\infty.$

For item (ii) the curve goes inside $U_\delta$ and outside $N$ 
infinite number of times. By item (c) we get 
\[
 \int_{t_k}^{s_{k}}|f||d\gamma|\geq 
\frac{1}{1+c_2^2}\int_{t_{k}}^{s_{k}}
|f|(1+|g|^2)|d\gamma|\geq\frac { c_1 } { 1+c_2^2 }.
\]

Then
\begin{equation}\label{eq1.01}
 \int_{t_k}^{s_{m+k}}|f||d\gamma|\geq 
\frac{1}{1+c_2^2}\sum_{j=0}^{j=m}\int_{t_{k+j}}^{s_{k+j}}
|f|(1+|g|^2)|d\gamma|\geq\frac { c_1 } { 1+c_2^2 } m ,
\end{equation}
implying $\tilde\ell(\gamma)=\infty.$

In item (iii) the set $V_\delta$ has finite connected components 
and all are precompact. Hence $c_3=\sup\text{dist}^M(\partial V_{\epsilon'},\ 
\partial  V_\delta)>0$ and $c_4=\sup_{z\in V_\delta\setminus 
V_{\delta'}}|G(z)|<\infty.$ In a similar way the proof of item (ii) we get 
\[\int_{t_k}^{s_{k+m}}|f||d\gamma|\geq\frac{c_3}{1+c_4^2}m.\] 
implying $\int_\gamma|f||dz|=\infty.$

Observe that $\phi^{-1}(y_0)$ is a sequence of points into $\D$ accumulating 
only at $\partial\D$.  If $\gamma$ is divergent in $\D$ we can not have 
$\lim_t\phi(\gamma(t))=y_0$ otherwise 
$\gamma(t)$ is in one connected component of $\phi^{-1}(V_\delta)$ for all 
$t\geq t_0>0$ and is not divergent. The unique way to have $y_0$ accumulation 
point of $\phi(\gamma)$ is if item (iii) happen. We already know that in this 
case $\tilde\ell(\gamma)=\infty$. This finish the proof of the 
lemma.

\end{proof} 

\begin{remark}\label{remark}
A complete minimal surfaces satisfying the conditions of lemma (\ref{lemma3}) 
does not exist.
 In fact we have proved that for all 
divergent 
curve $\gamma$ into $\D$ we have $\int_\gamma|f||dz|=\infty$. By a well 
know Osserman's lemma there is no holomorphic $f$ with finite number of zeros 
and $|f||dz|$ complete a 
metric in the unit disc (see lemma 8.5 of \cite{kn:Os}). Our $f$ has infinite 
number of zeros and do not exist for other reason.   
\end{remark}

\newpage


\end{document}